\documentclass[twoside,11pt,reqno]{amsart}
\usepackage{amsmath, amsthm,amssymb,amstext,amsfonts,amscd}
\usepackage{graphicx}
\usepackage{multirow}
\usepackage{lmodern}
\usepackage{latexsym}
\usepackage{wasysym} 
\usepackage{url}
\numberwithin{equation}{section}

\begin{document}

\title[{Some summation theorems for truncated Clausen series}]{Some summation theorems for truncated Clausen series and applications}

\author[M.I. Qureshi, Saima Jabee$^{*}$ and Dilshad Ahamad]{M.I. Qureshi, Saima Jabee$^{*}$ and Dilshad Ahamad}

\address{M.I. Qureshi: Department of Applied Sciences and Humanities,
 Faculty of Engineering and Technology,
 Jamia Millia Islamia (A Central University),
 New Delhi 110025, India.}
\email{miqureshi\_delhi@yahoo.co.in}

\address{Saima Jabee: Department of Applied Sciences and Humanities,
Faculty of Engineering and Technology,
Jamia Millia Islamia (A Central University),
New Delhi 110025, India. }
\email{saimajabee007@gmail.com}

\address{Dilshad Ahamad: Department of Applied Sciences and Humanities,
 Faculty of Engineering and Technology,
 Jamia Millia Islamia (A Central University),
 New Delhi 110025, India.}
\email{dlshdhmd4@gmail.com}



\thanks{*Corresponding author}

\begin{abstract}
The main aim of this paper is to derive some new summation theorems for terminating and truncated Clausen's hypergeometric series with unit argument, when one numerator parameter and one denominator parameter are negative integers. Further, using our truncated summation theorems, we obtain the Mellin transforms of the product of exponential function and Goursat's truncated hypergeometric function.
\end{abstract}
\begin{center}
\today
\end{center}
\maketitle
{
\section{Introduction}
In our investigations, we shall use the following standard notations:\\
$\mathbb{N}:=\{1,2,3,\dots\}$;  $\mathbb{N}_0:=\mathbb{N}\bigcup\{0\}$;  $\mathbb{Z}_0^-:=\mathbb{Z}^-\bigcup\{0\}=\{0,-1,-2,-3,\dots\}$.\\
The symbols $\mathbb{C}$, $\mathbb{R}$, $\mathbb{N}$, $\mathbb{Z}$, $\mathbb{R}^+$ and $\mathbb{R}^-$ denote the sets of complex numbers, real numbers, natural numbers, integers, positive and negative real numbers respectively.\\
The Pochhammer symbol $(\alpha)_{p}$ ~$(\alpha, p \in\mathbb{C})$ (\cite[p.22 eq(1), p.32 Q.N.(8) and Q.N.(9)]{Rainville}, see also \cite[p.23, eq(22) and eq(23)]{Srivastava3}), is defined by
\begin{equation}\label{f.eq(g1)}
(\alpha)_{p}:=\frac{\Gamma(\alpha+p)}{\Gamma(\alpha)}=
\begin{cases}
$1$ & ;(p=0; \alpha \in \mathbb {C}\setminus \{0\})\\
\alpha (\alpha+1)\ldots (\alpha+n-1) & ;(p=n \in \mathbb {N}; \alpha \in \mathbb {C}\setminus {\mathbb{Z}_0^-})\\
\frac{(-1)^{n}k!}{(k-n)!} & ;(\alpha=-k; p=n; n,k \in \mathbb{N}_0; {0}\leq{n}\leq{k})\\
$0$ & ;(\alpha=-k; p=n; n,k \in \mathbb{N}_0;{n}>{k})\\
\frac{(-1)^{n}}{(1-\alpha)_n} & ;(p=-n; n\in \mathbb{N}; \alpha \in \mathbb{C}\setminus \mathbb{Z}),
\end{cases}
\end{equation}
it being understood conventionally that $(0)_0=1$ and assumed tacitly that the Gamma quotient exists.\\
The generalized hypergeometric function ${_p}F_q$ (\cite[Art.44, pp.73-74]{Rainville}, see also \cite{Bailey}), is defined by
\begin{eqnarray}\label{f.eq(g2)}
{_p}F_{q}\left[\begin{array}{r} \alpha_1, \alpha_2, \dots, \alpha_p;\\
~\\ \beta_1, \beta_2, \dots, \beta_q;\end{array}\  z\right]={_p}F_{q}\left[\begin{array}{r} (\alpha_p);\\
~\\ (\beta_q);\end{array}\  z\right]=\sum_{n=0}^{\infty}\frac{\displaystyle\prod_{j=1}^{p}(\alpha_j)_n}{\displaystyle\prod_{j=1}^{q}(\beta_j)_n}\frac{z^n}{n!}.
\end{eqnarray}
By convention, a product over the empty set is unity.\\
$\big(p, q \in \mathbb{N}_0;~ p\leqq{q+1}~;~ p\leqq{q}~ \text{and}~ |z|<\infty ;\big.$
$~\big. p=q+1 ~\text{and}~ |z|<1;~ p=q+1, |z|=1~\text{and}~\Re(\omega)>0;~p=q+1, |z|=1, z\neq 1~\text{and}~ -1< \Re(\omega) \leq0\big)$,\\
where \[\omega:=\sum_{j=1}^{q}{\beta}_j-\sum_{j=1}^{p}{\alpha}_j, \] \[\big(\alpha_j\in \mathbb{C}~(j=1, 2,\dots,p ); \beta_j\in\mathbb{C}\setminus\mathbb{Z}_0^-(j=1, 2, \dots, q) \big),\]
where $\Re$ denotes the real part of complex number throughout the paper.\\
A finite series identity (reversal of the order of terms in finite summation) is given by
\begin{eqnarray}\label{eq(g14)}
\sum_{n=0}^{m}\Phi(n)=\sum_{n=0}^{m}\Phi(m-n);\quad{m\in\mathbb{N}_0}.
\end{eqnarray}
The truncated hypergeometric series is given by:
\begin{eqnarray}\label{f.eq(g16)}
&&\text{The sum of the first (m+1)-terms of infinite series } {_{p}}F_{q}\left[\begin{array}{r} (\alpha_p);\\
~\\ (\beta_q);\end{array}\  z\right]\nonumber\\
&&\qquad\qquad={_{p}}F_{q}\left[\begin{array}{r} (\alpha_p);\\
~\\ (\beta_q);\end{array}\  z\right]_{m}=\sum_{n=0}^{m}\frac{\displaystyle\prod_{j=1}^{p}(\alpha_j)_n}{\displaystyle\prod_{j=1}^{q}(\beta_j)_n}\frac{z^n}{n!}\nonumber\\
&&\qquad\qquad=\frac{[(\alpha_p)]_{m}z^{m}}{[(\beta_q)]_{m}m!}{_{q+2}}F_{p}\left[\begin{array}{r} -m, 1-(\beta_{q})-m,1;\\
~\\ 1-(\alpha_p)-m;\end{array}\  \frac{(-1)^{p+q+1}}{z}\right],
\end{eqnarray}
where $(\alpha_{p}), (\beta_{q}), 1-(\alpha_{p})-m, 1-(\beta_{q})-m\in\mathbb{C}\setminus\mathbb{Z}_0^-$; $m\in\mathbb{N}_0$, and
\begin{eqnarray}
[(\alpha_p)]_{m}=(\alpha_1)_{m}(\alpha_2)_{m}\dots(\alpha_p)_{m}=\prod_{i=1}^{p}(\alpha_i)_{m}=\prod_{i=1}^{p}\frac{\Gamma(\alpha_i+m)}{\Gamma(\alpha_i)},
\end{eqnarray}
with similar interpretation for others.\\
The terminating hypergeometric series (the hypergeometric polynomial) is given by
\begin{eqnarray}\label{eq(g15)}
{_{p+1}}F_{q}\left[\begin{array}{r} -m, (\alpha_p);\\
~\\ (\beta_q);\end{array}\ z\right]&=&\frac{[(\alpha_p)]_{m}(-z)^{m}}{[(\beta_q)]_{m}}{_{q+1}}F_{p}\left[\begin{array}{r} -m, 1-(\beta_{q})-m;\\
~\\ 1-(\alpha_p)-m;\end{array}\ \frac{(-1)^{p+q}}{z}\right],\nonumber\\
\end{eqnarray}
where $(\alpha_{p}), (\beta_{q}), 1-(\alpha_{p})-m, 1-(\beta_{q})-m\in\mathbb{C}\setminus\mathbb{Z}_0^-$ and $m\in\mathbb{N}_0$.\\

If $\ell> m$; $m,\ell\in\mathbb{N}; \alpha,\beta,\gamma\in\mathbb{C}\setminus\mathbb{Z}_0^{-}$, then series ${_3}F_2\left[\begin{array}{r} -m, \alpha,\beta;\\
~\\ -\ell,\gamma;\end{array}\ z\right]$ is an infinite series and is given by the following series representation (see for example \cite[p.41, eq.(3.1.26); p.42, eq.(3.2.6)]{Luke} and \cite[p.438, eq.(7.2.3.5)]{Prudnikov2})
\begin{eqnarray}
{_3}F_2\left[\begin{array}{r} -m, \alpha,\beta;\\
~\\ -\ell,\gamma;\end{array}\ z\right]&=&\sum_{r=0}^{m}\frac{(-m)_r(\alpha)_r(\beta)_rz^r}{(-\ell)_r(\gamma)_rr!}+\sum_{r=\ell+1}^{\infty}\frac{(-m)_r(\alpha)_r(\beta)_rz^r}{(-\ell)_r(\gamma)_rr!}\nonumber\\
&&={_3}F_2\left[\begin{array}{r} -m, \alpha,\beta;\\
~\\ -\ell,\gamma;\end{array}\ z\right]_m+\sum_{r=\ell+1}^{\infty}\frac{(-m)_r(\alpha)_r(\beta)_rz^r}{(-\ell)_r(\gamma)_rr!}.
\end{eqnarray}

In original notation, the higher order Goursat hypergeometric function is represented by double integral \cite[p. 286]{Goursat}. So we have
\begin{eqnarray}\label{eq(1)}
G\left(\begin{array}{r} \alpha,\beta;\\ \gamma, \delta;\end{array}z\right)&=&\frac{\Gamma{(\gamma)}\Gamma{(\delta)}}{\Gamma{(\alpha)}\Gamma{(\beta)}\Gamma{(\gamma-\alpha)}\Gamma{(\delta-\beta)}}\times\nonumber\\
&&\times\int_{0}^{1}\int_{0}^{1}u^{\alpha-1}v^{\beta-1}(1-u)^{\gamma-\alpha-1}(1-v)^{\delta-\beta-1}e^{zuv}dudv,
\end{eqnarray}
where $\Re{(\gamma)}>\Re{(\alpha)}>0$, $\Re{(\delta)}>\Re{(\beta)}>0$,\\
and
\begin{eqnarray*}
G\left(\begin{array}{r} \alpha,\beta;\\ \gamma, \delta;\end{array}z\right)&=&1+\sum_{n=1}^{\infty}\frac{(\alpha)_n(\beta)_n z^n}{(\gamma)_n(\delta)_n n!}\\
&=&{_2}F_{2}\left[\begin{array}{r} \alpha,\beta;\\ \gamma, \delta;\end{array}z\right],
\end{eqnarray*}
where $\gamma,\delta\in\mathbb{C}\setminus\mathbb{Z}_0^{-}$ and $|z|<\infty$.\\
It is also well known that, under certain conditions, the Goursat's function \cite[p. 286]{Goursat} ${_2}F_2(\alpha, \beta; \gamma, \delta; z)$ is defined by
\begin{eqnarray}\label{eq(2)}
{_2F_2} \left[\begin{array}{r} \alpha,\beta;\\ \gamma, \delta;\end{array}z\right]=\frac{\Gamma{(\delta)}}{\Gamma{(\alpha)}\Gamma{(\delta-\alpha)}}\int_{0}^{1}v^{\alpha-1}(1-v)^{\delta-\alpha-1}{_1F_1} \left[\begin{array}{r} \beta;\\ \gamma;\end{array}zv\right]dv,
\end{eqnarray}
where $\Re{(\delta)}>\Re{(\alpha)}>0$ and ${_1}F_{1}(\cdot)$ is Kummer's confluent hypergeometric function.\\
An integral transform that may be considered as the multiplicative form of the two-sided Laplace transform is known as Mellin transform, which is closely related to the Fourier transform, Laplace transform and other transforms. The Mellin transform is defined by
\begin{eqnarray}\label{eq(1.9)}
\mathcal{M}\{f(t);s\}=\int_{0}^{\infty}t^{s-1}f(t)dt=g(s),
\end{eqnarray}
where $s$ is a complex variable, above integral exists with suitable convergence conditions.\\
Until 1990, only few classical summation theorems for ${_2}F_1$ and ${_3}F_2$ were known. Subsequently, some progress has been made in generalizing these classical summation theorems (see \cite{Kim, Lavoie1, Lavoie2, Lavoie3, Miller, Rakha1, Rakha2}).\\
\section{Summation theorems for non-terminating, terminating and truncated clausen series }
In this section, we have verified the following terminating and truncated Clausen summation theorems by taking suitable values of parameters. So, without any loss of convergence, we can relax convergence conditions in some cases.\\
The classical Watson's summation theorem for non-terminating Clausen's hypergeometric series of unit argument \cite[p.16, section 3.3(1)]{Bailey} takes the form
\begin{eqnarray}\label{eq(2.1)}
{_3F_2} \left[\begin{array}{r} \alpha,\beta, \gamma;\\ \frac{1+\alpha+\beta}{2}, 2\gamma;\end{array}1\right] &=& \frac{\Gamma{\left(\frac{1}{2}\right)}\Gamma{\left(\gamma+\frac{1}{2}\right)}\Gamma{\left(\frac{1+\alpha+\beta}{2}\right)}\Gamma{\left(\gamma+\frac{1-\alpha-\beta}{2}\right)}}{\Gamma{\left(\frac{1+\alpha}{2}\right)}\Gamma{\left(\frac{1+\beta}{2}\right)}\Gamma{\left(\gamma+\frac{1-\alpha}{2}\right)}\Gamma{\left(\gamma+\frac{1-\beta}{2}\right)}},
\end{eqnarray}
provided $\Re(\gamma+\frac{1-\alpha-\beta}{2})>0; \frac{1+\alpha+\beta}{2}, \gamma, 2\gamma\in\mathbb{C}\setminus\mathbb{Z}_0^{-}$ and parameters are adjusted in such a way that the series on the left-hand side is well defined.\\
When $\alpha=-2m$ in equation \eqref{eq(2.1)}, we get a Watson's summation theorem for terminating hypergeometric series (containing (2m+1)-terms)
\begin{eqnarray}\label{eq(2.2)}
{_3F_2} \left[\begin{array}{r} -2m,\beta, \gamma;\\ \frac{1-2m+\beta}{2}, 2\gamma;\end{array}1\right] &=& \frac{\left(\frac{1}{2}\right)_m\left(\gamma+\frac{1-\beta}{2}\right)_m}{\left(\gamma+\frac{1}{2}\right)_m\left(\frac{1-\beta}{2}\right)_m},
\end{eqnarray}
where $\beta,\gamma, 2\gamma, \frac{1+\beta}{2}-m\in\mathbb{C}\setminus\mathbb{Z}_0^{-}$; $m\in\mathbb{N}$.\\
When $\alpha=-2m-1$ in equation \eqref{eq(2.1)}, we get another Watson's summation theorem for terminating hypergeometric series (containing-(2m+2) terms)
\begin{eqnarray}\label{eq(2.3)}
{_3F_2} \left[\begin{array}{r} -2m-1,\beta, \gamma;\\ \frac{\beta-2m}{2}, 2\gamma;\end{array}1\right] &=& 0,
\end{eqnarray}
where $\beta, \gamma, 2\gamma, \frac{\beta-2m}{2}\in\mathbb{C}\setminus\mathbb{Z}_0^{-}$; $m\in\mathbb{N}$.\\

We recall a Watson's summation theorem for truncated Clausen's series (containing (m+1)-terms) \cite[p.238, eq(2.2)]{Bailey1}
\begin{eqnarray}\label{eq(2.4)}
{_3F_2} \left[\begin{array}{r} -m,\alpha,\beta;\\ -2m,\frac{1+\alpha+\beta}{2};\end{array}1\right]_m &=& \frac{\left(\frac{1+\alpha}{2}\right)_m\left(\frac{1+\beta}{2}\right)_m}{\left(\frac{1}{2}\right)_m\left(\frac{1+\alpha+\beta}{2}\right)_m},
\end{eqnarray}
where $\alpha,\beta,\frac{1+\alpha+\beta}{2}\in\mathbb{C}\setminus\mathbb{Z}_0^{-}; m\in\mathbb{N}$.\\
On setting $\gamma=-m-k-\frac{1}{2}$ in equation \eqref{eq(2.2)}, we obtain Watson's summation theorem for truncated Clausen's series (containing-(2m+1) terms) is given by
\begin{eqnarray}\label{eq(2.5)}
{_3F_2} \left[\begin{array}{r} -2m,\beta,-m-k-\frac{1}{2};\\ -2m-2k-1, \frac{1+\beta}{2}-m;\end{array}1\right]_{2m}&=& \frac{\left(\frac{1}{2}\right)_m\left(\frac{2+\beta+2k}{2}\right)_m}{\left(\frac{1-\beta}{2}\right)_m\left(1+k\right)_m},
\end{eqnarray}
where $\beta,\frac{1+\beta}{2}-m\in\mathbb{C}\setminus\mathbb{Z}_0^{-}; m,k\in\mathbb{N}$.\\
On setting $\gamma=-m-k-\frac{1}{2}$ in equation \eqref{eq(2.3)}, we obtain  another Watson's summation theorem for truncated Clausen's series (containing-(2m+2) terms) is given by
\begin{eqnarray}\label{eq(2.6)}
{_3F_2} \left[\begin{array}{r} -2m-1,\beta,-m-k-\frac{1}{2};\\ -2m-2k-1, \frac{\beta}{2}-m;\end{array}1\right]_{2m+1}&=&0,
\end{eqnarray}
where $\beta,\frac{\beta}{2}-m\in\mathbb{C}\setminus\mathbb{Z}_0^{-}; m,k\in\mathbb{N}$.\\
The following summation theorem  for Clausen's non-terminating series due to Saalsch\"{u}tz's (\cite[p.21,section 3.8(2)]{Bailey}, \cite[p.534, Entry 12]{Prudnikov2}, see also \cite[p.73(2.4.4.4) and p.246(III.31)]{Slater}) is given by
\begin{eqnarray}\label{eq(2.7)}
&&{_3F_2} \left[\begin{array}{r} \alpha,\beta, \gamma+\delta-\alpha-
\beta-1;\\ \gamma, \delta;\end{array}1\right]\nonumber\\
&&=\frac{\Gamma{(\gamma)}\Gamma{(\delta)}\Gamma{(\gamma-\alpha-\beta)}\Gamma{(\delta-\alpha-\beta)}}{\Gamma{(\gamma-\alpha)}\Gamma{(\gamma-\beta)}\Gamma{(\delta-\alpha)}\Gamma{(\delta-\beta)}}+\frac{1}{(\alpha+\beta-\gamma)}\frac{\Gamma{(\gamma)}\Gamma{(\delta)}}{\Gamma{(\alpha)}\Gamma{(\beta)}\Gamma{(\gamma+\delta-\alpha-\beta)}}\times\nonumber\\
&&\times{_3F_2} \left[\begin{array}{r} \gamma-\alpha,\gamma-\beta, 1;\\ \gamma-\alpha-\beta+1, \gamma+\delta-\alpha-\beta;\end{array}1\right],
\end{eqnarray}
where $\Re{(\delta-\alpha-\beta)}>0$ and $\Re{(\gamma-\alpha-\beta)}>0$.\\
If we set $\delta=-m+1-\gamma+\alpha+\beta$, $m$ being positive integer, in the right-hand side of equation \eqref{eq(2.7)}, we obtain Saalsch\"{u}tz's summation theorem for Clausen's terminating series (\cite[p.9, section 2.2(1)]{Bailey}, see also \cite[p.87, Th 29]{Rainville})
\begin{eqnarray}\label{eq(2.8)}
{_3F_2} \left[\begin{array}{r} \alpha,\beta, -m;\\ \gamma, 1+\alpha+\beta-\gamma-m;\end{array}1\right] &=& \frac{\left(\gamma-\alpha\right)_m\left(\gamma-\beta\right)_m}{\left(\gamma\right)_m\left(\gamma-\alpha-\beta\right)_m},
\end{eqnarray}
where $ \alpha,\beta,\gamma, 1+\alpha+\beta-\gamma-m\in\mathbb{C}\setminus\mathbb{Z}_0^{-}; m\in\mathbb{N}$.\\
On setting $\gamma=-m-k$ in equation \eqref{eq(2.8)}, we get Saalsch\"{u}tz's summation theorem for truncated series
\begin{eqnarray}\label{eq(2.8c)}
{_3F_2} \left[\begin{array}{r} -m,\alpha,\beta;\\ -m-k, 1+\alpha+\beta+k;\end{array}1\right]_m &=& \frac{\left(1+\alpha+k\right)_m\left(1+\beta+k\right)_m}{\left(1+k\right)_m\left(1+\alpha+\beta+k\right)_m},
\end{eqnarray}
where $\alpha,\beta,1+\alpha+\beta+k\in\mathbb{C}\setminus\mathbb{Z}_0^{-}; m,k\in\mathbb{N}$.\\
Next we recall another Saalsch\"{u}tz's summation theorem for Clausen's terminating series (\cite[p.24]{Bailey1}, see also \cite[p.87, Theorem 30]{Rainville})
\begin{eqnarray}\label{eq(2.8a)}
{_3F_2} \left[\begin{array}{r} -m, \alpha+m, 1+\alpha-\beta-\gamma;\\ 1+\alpha-\beta, 1+\alpha-\gamma;\end{array}1\right] &=& \frac{\left(\beta\right)_m\left(\gamma\right)_m}{\left(1+\alpha-\beta\right)_m\left(1+\alpha-\gamma\right)_m},
\end{eqnarray}
where $\alpha+m, 1+\alpha-\beta-\gamma,1+\alpha-\beta,1+\alpha-\gamma\in\mathbb{C}\setminus\mathbb{Z}_0^{-}; m\in\mathbb{N}$.\\
If we set $1+\alpha-\beta=-m-k$ in equation \eqref{eq(2.8a)}, we obtain the following Saalsch\"{u}tz's summation theorem for truncated series
\begin{eqnarray}\label{eq(2.8b)}
{_3F_2} \left[\begin{array}{r} -m, \beta-k-1, -m-k-\gamma;\\ -m-k, \beta-\gamma-m-k;\end{array}1\right]_m &=& \frac{\left(\beta\right)_m\left(\gamma\right)_m}{\left(1+k\right)_m\left(1+k+\gamma-\beta\right)_m},
\end{eqnarray}
where $\beta-k-1,-m-k-\gamma,\beta-\gamma-m-k\in\mathbb{C}\setminus\mathbb{Z}_0; m,k\in\mathbb{N}$.\\
Next we recall Whipple's summation theorem for non-terminating Clausen's series \cite[p.16, section 3.4(1)]{Bailey}
\begin{eqnarray}\label{eq(2.9)}
&&{_3F_2} \left[\begin{array}{r} \alpha,1- \alpha, \beta;\\ \gamma, 2\beta-\gamma+1;\end{array}1\right]\nonumber\\
&&=\frac{\pi \Gamma{(\gamma)}\Gamma{(2\beta-\gamma+1)}}
{2^{2\beta-1}\Gamma{\left(\frac{\alpha+2\beta-\gamma+1}{2}\right)}\Gamma{\left(\frac{\alpha+\gamma}{2}\right)}\Gamma{\left(\frac{2-\alpha+2\beta-\gamma}{2}\right)}\Gamma{\left(\frac{1-\alpha+\gamma}{2}\right)}},
\end{eqnarray}
where $\Re{(\beta)}>0,\gamma,2\beta-\gamma+1\in\mathbb{C}\setminus\mathbb{Z}_0^{-}$.\\
On setting $\alpha=-2m$ in equation \eqref{eq(2.9)}, we get Whipple's summation theorem for terminating series
\begin{eqnarray}\label{eq(2.10)}
{_3F_2} \left[\begin{array}{r} -2m,1+2m, \beta;\\ \gamma, 1+2\beta-\gamma;\end{array}1\right]=\frac{\left(\frac{2-\gamma}{2}\right)_m\left(\frac{1-2\beta+\gamma}{2}\right)_m}{\left(\frac{1+\gamma}{2}\right)_m\left(\frac{2+2\beta-\gamma}{2}\right)_m},
\end{eqnarray}
where $\beta,\gamma,1+2\beta-\gamma\in\mathbb{C}\setminus\mathbb{Z}_0^{-}$; $m\in\mathbb{N}$.\\
On setting $\alpha=-2m-1$ in equation \eqref{eq(2.9)}, we get another Whipple's summation theorem for terminating series
\begin{eqnarray}\label{eq(2.11)}
{_3F_2} \left[\begin{array}{r} -2m-1,2+2m, \beta;\\ \gamma, 1+2\beta-\gamma;\end{array}1\right]=\frac{(\gamma-1)(2\beta-\gamma)\left(\frac{3-\gamma}{2}\right)_m\left(\frac{2-2\beta+\gamma}{2}\right)_m}{(\gamma)(1+2\beta-\gamma)\left(\frac{2+\gamma}{2}\right)_m\left(\frac{3+2\beta-\gamma}{2}\right)_m},
\end{eqnarray}
where $\beta,\gamma,1+2\beta-\gamma\in\mathbb{C}\setminus\mathbb{Z}_0^{-}$; $m\in\mathbb{N}$.\\
Another Whipple's summation theorem for Clausen's terminating series is given by (\cite[p. 157, eq(3.1)]{Qureshi4}, see also \cite[p. 190, eq(2)]{Dzhrbashyan} and \cite[p. 238, eq(3.1)]{Bailey1})
\begin{eqnarray}\label{eq(2.13)}
&&{_3F_2} \left[\begin{array}{r} -m,\alpha, 1-\alpha;\\ \gamma, 1-\gamma-2m;\end{array}1\right]=\frac{\left(\frac{\gamma+\alpha}{2}\right)_{m} \left(\frac{\gamma-\alpha+1}{2}\right)_{m}}{\left(\frac{\gamma}{2}\right)_{m} \left(\frac{\gamma+1}{2}\right)_{m}},
\end{eqnarray}
where $\alpha,1-\alpha, \gamma, 1-\gamma-2m, \frac{\alpha+\gamma+2m}{2},\frac{1-\alpha+\gamma+2m}{2}\in\,\mathbb{C}\setminus\mathbb{Z}_0^-$; $m\in\mathbb{N}$.\\
If we set $\gamma=-2m-k$ in equation \eqref{eq(2.13)}, we get Whipple summation theorem for truncated series containing (m+1)-terms
\begin{eqnarray}\label{eq(2.14a)}
{_3F_2} \left[\begin{array}{r} -m,\alpha, 1-\alpha;\\ -2m-k, 1+k;\end{array}1\right]_m=\frac{\left(\frac{2-\alpha+k}{2}\right)_{m} \left(\frac{1+\alpha+k}{2}\right)_{m}}{\left(\frac{2+k}{2}\right)_{m} \left(\frac{1+k}{2}\right)_{m}},
\end{eqnarray}
where $\alpha,1-\alpha\in\mathbb{C}\setminus\mathbb{Z}_0^{-}; m,k\in\mathbb{N}$.\\
On setting $\gamma=-2m-2k$ in equation \eqref{eq(2.10)}, we get Whipple summation theorem for truncated series containing (2m+1)-terms
\begin{eqnarray}\label{eq(2.15)}
&&{_3F_2} \left[\begin{array}{r} -2m,1+2m,\beta;\\ -2m-2k, 2\beta+2m+2k+1;\end{array}1\right]_{2m}=\frac{(1+2\beta+2k)_{2m}\left(1+k\right)_{2m}}{(1+2k)_{2m}\left(1+\beta+k\right)_{2m}},
\end{eqnarray}
where $\beta,2\beta+2m+2k+1\in\mathbb{C}\setminus\mathbb{Z}_0^{-}; m,k\in\mathbb{N}$.\\
On setting $\gamma=-2m-2k-1$ in equation \eqref{eq(2.10)}, we get Whipple summation theorem for truncated series containing (2m+1)-terms
\begin{eqnarray}\label{eq(2.15a)}
&&{_3F_2} \left[\begin{array}{r} -2m,1+2m,\beta;\\ -2m-2k-1, 2\beta+2+2m+2k;\end{array}1\right]_{2m}=\frac{(2+2\beta+2k)_{2m}\left(\frac{3+2k}{2}\right)_{2m}}{(2+2k)_{2m}\left(\frac{3+2\beta+2k}{2}\right)_{2m}},\nonumber\\
\end{eqnarray}
where $\beta,2\beta+2m+2k+2\in\mathbb{C}\setminus\mathbb{Z}_0^{-}; m,k\in\mathbb{N}$.\\
If we set $\gamma=-2m-2k-1$ in equation \eqref{eq(2.11)}, we get Whipple summation theorem for truncated series containing (2m+2)-terms
\begin{eqnarray}\label{eq(2.16)}
&&{_3F_2} \left[\begin{array}{r} -2m-1,2+2m,\beta;\\ -2m-2k-1, 2\beta+2m+2k+2;\end{array}1\right]_{2m+1}\nonumber\\
&&=\frac{(k+1)(2\beta+2m+2k+1)(2\beta+2k+1)_{2m}\left(2+k\right)_{2m}}{(2m+2k+1)(\beta+k+1)(2k+1)_{2m}\left(2+\beta+k\right)_{2m}},
\end{eqnarray}
where $\beta,2\beta+2m+2k+2\in\mathbb{C}\setminus\mathbb{Z}_0^{-}; m,k\in\mathbb{N}$.\\
If we set $\gamma=-2m-2k-2$ in equation \eqref{eq(2.11)}, we get Whipple summation theorem for truncated series containing (2m+2)-terms
\begin{eqnarray}\label{eq(2.16a)}
&&{_3F_2} \left[\begin{array}{r} -2m-1,2+2m,\beta;\\ -2m-2k-2, 2\beta+2m+2k+3;\end{array}1\right]_{2m+1}\nonumber\\
&&=\frac{(2k+3)(\beta+m+k+1)(2\beta+2k+2)_{2m}\left(\frac{5+2k}{2}\right)_{2m}}{(m+k+1)(2\beta+2k+3)(2k+2)_{2m}\left(\frac{5+2\beta+2k}{2}\right)_{2m}},
\end{eqnarray}
where $\beta,2\beta+2m+2k+3\in\mathbb{C}\setminus\mathbb{Z}_0^{-}; m,k\in\mathbb{N}$.\\
The classical Dixon's summation theorem for Clausen's non-terminating series \cite[p.13, section 3.1(1)]{Bailey} is given by
\begin{eqnarray}\label{eq(2.18)}
&&{_3F_2} \left[\begin{array}{r} \alpha,\beta, \gamma;\\ 1+\alpha-\beta, 1+\alpha-\gamma;\end{array}1\right]\nonumber\\
&&=\frac{\Gamma{\left(1+\frac{\alpha}{2}\right)}\Gamma{(1+\alpha-\beta)}\Gamma{(1+\alpha-\gamma)}\Gamma{\left(1+\frac{\alpha}{2}-\beta-\gamma\right)}}{\Gamma{\left(1+\alpha\right)}\Gamma{\left(1+\frac{\alpha}{2}-\beta\right)}\Gamma{\left(1+\frac{\alpha}{2}-\gamma\right)}\Gamma{\left(1+\alpha-\beta-\gamma\right)}},
\end{eqnarray}
where $\Re{(\alpha-2\beta-2\gamma)}>-2; \alpha,\beta,\gamma\in\mathbb{C};1+\alpha-\beta,1+\alpha-\gamma,1+\frac{\alpha}{2},1+\frac{\alpha}{2}-\beta-\gamma\in\mathbb{C}\setminus\mathbb{Z}_0^-$.\\
Equation \eqref{eq(2.18)} can be written as
\begin{eqnarray}\label{eq(2.19)}
&&{_3F_2} \left[\begin{array}{r} \alpha,\beta, \gamma;\\ 1+\alpha-\beta, 1+\alpha-\gamma;\end{array}1\right]\nonumber\\
&&=\frac{\cos\left(\frac{\pi \alpha}{2}\right)\Gamma{(1-\alpha)}\Gamma{(1+\alpha-\beta)}\Gamma{(1+\alpha-\gamma)}\Gamma{\left(1+\frac{\alpha}{2}-\beta-\gamma\right)}}{\Gamma{\left(1-\frac{\alpha}{2}\right)}\Gamma{\left(1+\frac{\alpha}{2}-\beta\right)}\Gamma{\left(1+\frac{\alpha}{2}-\gamma\right)}\Gamma{\left(1+\alpha-\beta-\gamma\right)}}\\
&&=\frac{\cos\left(\frac{\pi \alpha}{2}\right)\Gamma{\left(\beta-\frac{\alpha}{2}\right)}\Gamma{(\gamma-\frac{\alpha}{2})}\Gamma{(1-\alpha)}\Gamma{\left(\beta+\gamma-\alpha\right)}}{\Gamma{\left(\beta-\alpha\right)}\Gamma{(\gamma-\alpha)}\Gamma{\left(1-\frac{\alpha}{2}\right)}\Gamma{\left(\beta+\gamma-\frac{\alpha}{2}\right)}}\times\nonumber\\
&&\times\frac{\sin\{\pi \left(\beta-\frac{\alpha}{2}\right)\}\sin\{\pi \left(\gamma-\frac{\alpha}{2}\right)\}\sin\{\pi \left(\beta+\gamma-\alpha\right)\}}{\sin\{\pi \left(\beta-\alpha\right)\}\sin\{\pi \left(\gamma-\alpha\right)\}\sin\{\pi \left(\beta+\gamma-\frac{\alpha}{2}\right)\}}.
\end{eqnarray}
On setting $\alpha=-2m$ in equation \eqref{eq(2.18)}, we obtain Dixon's summation theorem for terminating series
\begin{eqnarray}\label{eq(2.20)}
{_3F_2} \left[\begin{array}{r} -2m,\beta, \gamma;\\ 1-2m-\beta, 1-2m-\gamma;\end{array}1\right]=\frac{(\beta)_{m}(\gamma)_{m}2^{2m}\left(\frac{1}{2}\right)_m(\beta+\gamma)_{2m}}{(\beta)_{2m}(\gamma)_{2m}(\beta+\gamma)_{m}},
\end{eqnarray}
where $\beta, \gamma\in \mathbb{C}\setminus\mathbb{Z}; m\in\mathbb{N}$.\\
On setting $\alpha=-2m-1$ in equation \eqref{eq(2.18)}, we obtain another Dixon's summation theorem for terminating series
\begin{eqnarray}\label{eq(2.21)}
&&{_3F_2} \left[\begin{array}{r} -2m-1,\beta, \gamma;\\ -2m-\beta, -2m-\gamma;\end{array}1\right]=0,
\end{eqnarray}
where $\beta, \gamma\in \mathbb{C}\setminus\mathbb{Z}; m\in\mathbb{N}$.\\
On setting $\beta=1+k$ in equation \eqref{eq(2.20)}, we obtain Dixon's summation theorem for truncated series
\begin{eqnarray}\label{eq(2.22)}
{_3F_2} \left[\begin{array}{r} -2m,1+k, \gamma;\\ -2m-k, 1-2m-\gamma;\end{array}1\right]_{2m}=\frac{(1+k)_{m}(\gamma)_{m}2^{2m}\left(\frac{1}{2}\right)_m(1+k+\gamma)_{2m}}{(1+k)_{2m}(\gamma)_{2m}(1+k+\gamma)_{m}},
\end{eqnarray}
where $\gamma,1-2m-\gamma\in\mathbb{C}\setminus\mathbb{Z}_0^{-}; m,k\in\mathbb{N}$.\\
On setting $\gamma=1+k$ in equation \eqref{eq(2.22)}, we obtain another Dixon's summation theorem for truncated series
\begin{eqnarray}\label{eq(2.23)}
{_3F_2} \left[\begin{array}{r} -2m,1+k, 1+k;\\ -2m-k, -2m-k;\end{array}1\right]_{2m}=\frac{(1+k)_{m}(1+k)_{m}2^{2m}\left(\frac{1}{2}\right)_m(2+2k)_{2m}}{(1+k)_{2m}(1+k)_{2m}(2+2k)_{m}},
\end{eqnarray}
where $m,k\in\mathbb{N}$.\\
On setting $\beta=1+k$ in equation \eqref{eq(2.21)}, we obtain Dixon's summation theorem for truncated series
\begin{eqnarray}\label{eq(2.24)}
&&{_3F_2} \left[\begin{array}{r} -2m-1,1+k, \gamma;\\ -2m-1-k, -2m-\gamma;\end{array}1\right]_{2m+1}=0,
\end{eqnarray}
where $\gamma,-2m-\gamma\in\mathbb{C}\setminus\mathbb{Z}_0^{-}; m,k\in\mathbb{N}$.\\
On setting $\gamma=1+k$ in equation \eqref{eq(2.24)}, we obtain another Dixon's summation theorem for truncated series
\begin{eqnarray}\label{eq(2.25)}
&&{_3F_2} \left[\begin{array}{r} -2m-1,1+k, 1+k;\\ -2m-1-k, -2m-1-k;\end{array}1\right]_{2m+1}=0,
\end{eqnarray}
where $m,k\in\mathbb{N}$.\\
On setting $\beta=-2m$ in equation \eqref{eq(2.18)}, we obtain Dixon's summation theorem for terminating series
\begin{eqnarray}\label{eq(2.20a)}
{_3F_2} \left[\begin{array}{r} -2m,\alpha, \gamma;\\ 1+\alpha+2m, 1+\alpha-\gamma;\end{array}1\right]=\frac{(1+\alpha)_{2m}\left(1+\frac{\alpha}{2}-\gamma\right)_{2m}}{\left(1+\frac{\alpha}{2}\right)_{2m}(1+\alpha-\gamma)_{2m}},
\end{eqnarray}
where $\alpha, \gamma\in \mathbb{C}\setminus\mathbb{Z}; m\in\mathbb{N}$.\\
On setting $\beta=-2m-1$ in equation \eqref{eq(2.18)}, we obtain another Dixon's summation theorem for terminating series
\begin{eqnarray}\label{eq(2.20b)}
{_3F_2} \left[\begin{array}{r} -2m-1,\alpha, \gamma;\\ 2+\alpha+2m, 1+\alpha-\gamma;\end{array}1\right]=\frac{(1+\alpha)(2+\alpha-2\gamma)(2+\alpha)_{2m}\left(2+\frac{\alpha}{2}-\gamma\right)_{2m}}{(2+\alpha)(1+\alpha-\gamma)\left(2+\frac{\alpha}{2}\right)_{2m}(2+\alpha-\gamma)_{2m}},
\end{eqnarray}
where $\alpha, \gamma\in \mathbb{C}\setminus\mathbb{Z}; m\in\mathbb{N}$.\\
On setting $\gamma=1+\alpha+2m+k$ in equation \eqref{eq(2.20a)}, we obtain Dixon's summation theorem for truncated series
\begin{eqnarray}\label{eq(2.20c)}
{_3F_2} \left[\begin{array}{r} -2m,\alpha, 1+\alpha+2m+k;\\ -2m-k, 1+\alpha+2m;\end{array}1\right]_{2m}=\frac{(1+\alpha)_{2m}\left(1+\frac{\alpha}{2}+k\right)_{2m}}{\left(1+\frac{\alpha}{2}\right)_{2m}(1+k)_{2m}},
\end{eqnarray}
where $\alpha,1+\alpha+2m+k,1+\alpha+2m\in\mathbb{C}\setminus\mathbb{Z}_0^{-}; m,k\in\mathbb{N}$.\\

On setting $\gamma=2+\alpha+2m+k$ in equation \eqref{eq(2.20b)}, we obtain another Dixon's summation theorem for truncated series
\begin{eqnarray}\label{eq(2.20d)}
&&{_3F_2} \left[\begin{array}{r} -2m-1,\alpha, 2+\alpha+2m+k;\\ -2m-k-1, 2+\alpha+2m;\end{array}1\right]_{2m+1}\nonumber\\
&&=\frac{(1+\alpha)(2+2k+\alpha+4m)(2+\alpha)_{2m}\left(1+\frac{\alpha}{2}+k\right)_{2m}}{(2+\alpha)(1+2m+k)\left(2+\frac{\alpha}{2}\right)_{2m}(1+k)_{2m}},
\end{eqnarray}
where $\alpha,2+\alpha+2m+k,2+\alpha+2m\in\mathbb{C}\setminus\mathbb{Z}_0^{-}; m,k\in\mathbb{N}$.\\

Also, on setting $\gamma=-m$ in equation \eqref{eq(2.19)}, we obtain Dixon's theorem for Clausen's terminating series
\begin{eqnarray}\label{eq(2.26)}
&&{_3F_2} \left[\begin{array}{r} \alpha, \beta, -m;\\ 1+\alpha-\beta, 1+\alpha+m;\end{array}1\right]\nonumber\\
&&=\frac{\cos\left(\frac{\pi \alpha}{2}\right)\Gamma{(1-\alpha)}\Gamma{(1+\alpha-\beta)}\Gamma{(1+\alpha+m)}\Gamma{\left(1+\frac{\alpha}{2}-\beta+m\right)}}{\Gamma{\left(1-\frac{\alpha}{2}\right)}\Gamma{\left(1+\frac{\alpha}{2}-\beta\right)}\Gamma{\left(1+\frac{\alpha}{2}+m\right)}\Gamma{\left(1+\alpha-\beta+m\right)}},
\end{eqnarray}
where $\alpha,\beta,1+\alpha-\beta,1+\alpha+m\in\mathbb{C}\setminus\mathbb{Z}_0^-; m\in \mathbb{N}$.\\
On setting $\beta=-m$ in equation \eqref{eq(2.18)}, we get Dixon's summation theorem for Clausen's terminating series
\begin{eqnarray}\label{eq(2.27)}
{_3F_2} \left[\begin{array}{r} -m, \alpha, \gamma;\\ 1+\alpha+m, 1+\alpha-\gamma;\end{array}1\right]=\frac{(1+\alpha)_m\left(1+\frac{\alpha}{2}-\gamma\right)_m}{\left(1+\frac{\alpha}{2}\right)_m(1+\alpha-\gamma)_m},
\end{eqnarray}
where $\alpha,\gamma,1+\alpha+m,1+\alpha-\gamma\in\mathbb{C}\setminus\mathbb{Z}_0^-; m\in \mathbb{N}$.\\
In section 3, we discuss the applications of some summation theorems for truncated Clausen hypergeometric series in Mellin transforms of the product of exponential function and truncated Goursat hypergeometric function.
\section{Applications in Mellin transforms}
In this section, we obtain Mellin transforms of the product of exponential function and truncated Goursat's function ${_2F_2}(\cdot)$ (when one numerator and one denominator parameters are negative integers),
\begin{eqnarray}\label{eq(7.1)}
\mathcal{M}\left\{e^{-\mu t}{_2F_2} \left[\begin{array}{r} -m, a;\\ -m-\ell, b;\end{array}\lambda t\right]_m;s\right\}&=&\int_{0}^{\infty}t^{s-1}e^{-\mu t}{_2F_2} \left[\begin{array}{r} -m, a;\\ -m-\ell, b;\end{array}\lambda t\right]_mdt\nonumber\\
&&=\frac{\Gamma(s)}{\mu^s}{_3F_2} \left[\begin{array}{r} -m, a,s;\\ -m-\ell, b;\end{array}\frac{\lambda}{\mu}\right]_m,
\end{eqnarray}
where $ \Re{(s)}>0$; $\Re{(\mu)}>0$ and $m,\ell\in\mathbb{N}$.\\
We derive some new results for Mellin transform as applications of summation theorems discussed in previous section.\\
\textbf{Case I.} On setting $\ell=m, a=\alpha, b=\frac{1+\alpha+\beta}{2}, \lambda=\mu$ and $s=\beta$ in equation \eqref{eq(7.1)} and using Watson's truncated summation theorem \eqref{eq(2.4)}, we obtain
\begin{eqnarray}\label{eq(7.2)}
\mathcal{M}\left\{e^{-\mu t}{_2F_2} \left[\begin{array}{r} -m,\alpha;\\ -2m,\frac{1+\alpha+\beta}{2};\end{array}\mu t\right]_m;\beta\right\}=\frac{\Gamma(\beta)}{\mu^{\beta}}\frac{\left(\frac{1+\alpha}{2}\right)_m\left(\frac{1+\beta}{2}\right)_m}{\left(\frac{1}{2}\right)_m\left(\frac{1+\alpha+\beta}{2}\right)_m},
\end{eqnarray}
where $\alpha,\frac{1+\alpha+\beta}{2}\in\mathbb{C}\setminus\mathbb{Z}_0^{-}$; $m\in\mathbb{N}$; $\Re(\beta)>0,\Re(\mu)>0$.\\

\textbf{Case II.} Replacing $m$ by $2m$ and after that setting $\ell=2k+1,a=-m-k-\frac{1}{2}, b=\frac{1+\beta}{2}-m, \lambda=\mu$ and $s=\beta$ in equation \eqref{eq(7.1)} and using Watson's truncated summation theorem \eqref{eq(2.5)}, we obtain
\begin{eqnarray}\label{eq(7.14)}
\mathcal{M}\left\{e^{-\mu t}{_2F_2} \left[\begin{array}{r} -2m,-m-k-\frac{1}{2};\\ -2m-2k-1, \frac{1+\beta}{2}-m;\end{array}\mu t\right]_{2m};\beta\right\}=\frac{\Gamma(\beta)}{\mu^{\beta}}\frac{\left(\frac{1}{2}\right)_m\left(\frac{2+\beta+2k}{2}\right)_m}{\left(\frac{1-\beta}{2}\right)_m\left(1+k\right)_m},
\end{eqnarray}
where $\left(\frac{1+\beta}{2}\right)-m\in\mathbb{C}\setminus\mathbb{Z}_0^{-}; m,k\in\mathbb{N};$ $\Re(\beta)>0, \Re(\mu)>0$.\\

\textbf{Case III.} Replacing $m$ by $2m+1$ and after that setting $\ell=2k,a=-m-k-\frac{1}{2}, b=\frac{\beta}{2}-m, \lambda=\mu$ and $s=\beta$ in equation \eqref{eq(7.1)} and using Watson's truncated summation theorem \eqref{eq(2.6)}, we obtain
\begin{eqnarray}\label{eq(7.15)}
\mathcal{M}\left\{e^{-\mu t}{_2F_2} \left[\begin{array}{r} -2m-1,-m-k-\frac{1}{2};\\ -2m-2k-1, \frac{\beta}{2}-m;\end{array}\mu t\right]_{2m+1};\beta\right\}=0,
\end{eqnarray}
where $\frac{\beta}{2}-m\in\mathbb{C}\setminus\mathbb{Z}_0^{-}; m,k\in\mathbb{N};$ $\Re(\beta)>0, \Re(\mu)>0$.\\

\textbf{Case IV.} On setting $\ell=k,a=\alpha, b=1+\alpha+\beta+k, \lambda=\mu$ and $s=\beta$ in equation \eqref{eq(7.1)} and using Saalsch\"{u}tz's truncated summation theorem \eqref{eq(2.8c)}, we obtain
\begin{eqnarray}\label{eq(7.3)}
\mathcal{M}\left\{e^{-\mu t}{_2F_2} \left[\begin{array}{r} -m,\alpha;\\ -m-k, 1+\alpha+\beta+k;\end{array}\mu t\right]_m;\beta\right\}= \frac{\Gamma(\beta)}{\mu^{\beta}}\frac{\left(1+\alpha+k\right)_m\left(1+\beta+k\right)_m}{\left(1+k\right)_m\left(1+\alpha+\beta+k\right)_m},\nonumber\\
\end{eqnarray}
where $\alpha,1+\alpha+\beta+k\in\mathbb{C}\setminus\mathbb{Z}_0^{-}; m,k\in\mathbb{N};$ $\Re(\beta)>0, \Re(\mu)>0$.\\

\textbf{Case V.} On setting $\ell=k,a=\beta-k-1, b=\beta-\gamma-m-k, \lambda=\mu$ and $s=-m-k-\gamma$ in equation \eqref{eq(7.1)} and using Saalsch\"{u}tz's truncated summation theorem \eqref{eq(2.8b)}, we obtain
\begin{eqnarray}\label{eq(g7.4)}
&&\mathcal{M}\left\{e^{-\mu t}{_2F_2} \left[\begin{array}{r} -m,\beta-k-1;\\ -m-k, \beta-\gamma-m-k;\end{array}\mu t\right]_m;-m-k-\gamma\right\}\nonumber\\
&&\qquad=\frac{\Gamma(-m-k-\gamma)}{\mu^{-m-k-\gamma}}\frac{\left(\beta\right)_m\left(\gamma\right)_m}{\left(1+k\right)_m\left(1+k+\gamma-\beta\right)_m},
\end{eqnarray}
where $\beta-k-1,\beta-\gamma-m-k\in\mathbb{C}\setminus\mathbb{Z}_0^{-}; m,k\in\mathbb{N};$ $\Re(-m-k-\gamma)>0, \Re(\mu)>0$.\\

\textbf{Case VI.} On setting $\ell=m+k,a=1-\alpha, b=1+k, \lambda=\mu$ and $s=\alpha$ in equation \eqref{eq(7.1)} and using Whipple's truncated summation theorem \eqref{eq(2.14a)}, we obtain
\begin{eqnarray}\label{eq(7.5)}
&&\mathcal{M}\left\{e^{-\mu t}{_2F_2} \left[\begin{array}{r} -m,1-\alpha;\\ -2m-k, 1+k;\end{array}\mu t\right]_m;\alpha\right\}=\frac{\Gamma(\alpha)}{\mu^{\alpha}}\frac{\left(\frac{2-\alpha+k}{2}\right)_{m} \left(\frac{1+\alpha+k}{2}\right)_{m}}{\left(\frac{2+k}{2}\right)_{m} \left(\frac{1+k}{2}\right)_{m}},
\end{eqnarray}
where $1-\alpha\in\mathbb{C}\setminus\mathbb{Z}_0^{-}; m,k\in\mathbb{N}$ and $\Re(\alpha)>0,\Re(\mu)>0$.\\

\textbf{Case VII.} Replacing $m$ by $2m$ and after that setting $\ell=2k,a=1+2m, b=1+2m+2k+2\beta, \lambda=\mu$ and $s=\beta$ in equation \eqref{eq(7.1)} and using Whipple's truncated summation theorem \eqref{eq(2.15)}, we obtain
\begin{eqnarray}\label{eq(7.6)}
&&\mathcal{M}\left\{e^{-\mu t}{_2F_2} \left[\begin{array}{r} -2m,1+2m;\\ -2m-2k, 2\beta+1+2m+2k;\end{array}\mu t\right]_{2m};\beta\right\}\nonumber\\
&&\qquad=\frac{\Gamma(\beta)}{\mu^{\beta}}\frac{(1+2\beta+2k)_{2m}\left(1+k\right)_{2m}}{(1+2k)_{2m}\left(1+\beta+k\right)_{2m}},
\end{eqnarray}
where $1+2m+2k+2\beta\in\mathbb{C}\setminus\mathbb{Z}_0^{-}; m,k\in\mathbb{N}$ and $\Re(\beta)>0,\Re(\mu)>0$.\\

\textbf{Case VIII.} Replacing $m$ by $2m$ and after that setting $\ell=2k+1,a=1+2m, b=2+2m+2k+2\beta, \lambda=\mu$ and $s=\beta$ in equation \eqref{eq(7.1)} and using Whipple's truncated summation theorem \eqref{eq(2.15a)}, we obtain
\begin{eqnarray}\label{eq(7.7)}
&&\mathcal{M}\left\{e^{-\mu t}{_2F_2} \left[\begin{array}{r} -2m,1+2m;\\ -2m-2k-1, 2\beta+2+2m+2k;\end{array}\mu t\right]_{2m};\beta\right\}\nonumber\\
&&\qquad=\frac{\Gamma(\beta)}{\mu^{\beta}}\frac{(2+2\beta+2k)_{2m}\left(\frac{3+2k}{2}\right)_{2m}}{(2+2k)_{2m}\left(\frac{3+2\beta+2k}{2}\right)_{2m}},
\end{eqnarray}
where $2+2m+2k+2\beta\in\mathbb{C}\setminus\mathbb{Z}_0^{-}; m,k\in\mathbb{N}$ and $\Re(\beta)>0,\Re(\mu)>0$.\\

\textbf{Case IX.} Replacing $m$ by $2m+1$ and after that setting $\ell=2k,a=2+2m, b=2+2m+2k+2\beta, \lambda=\mu$ and $s=\beta$ in equation \eqref{eq(7.1)} and using Whipple's truncated summation theorem \eqref{eq(2.16)}, we obtain
\begin{eqnarray}\label{eq(7.8)}
&&\mathcal{M}\left\{e^{-\mu t}{_2F_2} \left[\begin{array}{r} -2m-1,2+2m;\\ -2m-2k-1, 2\beta+2m+2k+2;\end{array}\mu t\right]_{2m+1};\beta\right\}\nonumber\\
&&=\frac{\Gamma(\beta)}{\mu^{\beta}}\frac{(k+1)(2\beta+2m+2k+1)(2\beta+2k+1)_{2m}\left(2+k\right)_{2m}}{(2m+2k+1)(\beta+k+1)(2k+1)_{2m}\left(2+\beta+k\right)_{2m}},
\end{eqnarray}
where $2+2m+2k+2\beta\in\mathbb{C}\setminus\mathbb{Z}_0^{-}; m,k\in\mathbb{N}$ and $\Re(\beta)>0,\Re(\mu)>0$.\\

\textbf{Case X.} Replacing $m$ by $2m+1$ and after that setting $\ell=2k+1,a=2+2m, b=3+2m+2k+2\beta, \lambda=\mu$ and $s=\beta$ in equation \eqref{eq(7.1)} and using Whipple's truncated summation theorem \eqref{eq(2.16a)}, we obtain
\begin{eqnarray}\label{eq(7.9)}
&&\mathcal{M}\left\{e^{-\mu t}{_2F_2} \left[\begin{array}{r} -2m-1,2+2m;\\ -2m-2k-2, 2\beta+2m+2k+3;\end{array}\mu t\right]_{2m+1};\beta\right\}\nonumber\\
&&=\frac{\Gamma(\beta)}{\mu^{\beta}}\frac{(2k+3)(\beta+m+k+1)(2\beta+2k+2)_{2m}\left(\frac{5+2k}{2}\right)_{2m}}{(m+k+1)(2\beta+2k+3)(2k+2)_{2m}\left(\frac{5+2\beta+2k}{2}\right)_{2m}},
\end{eqnarray}
where $3+2m+2k+2\beta\in\mathbb{C}\setminus\mathbb{Z}_0^{-}; m,k\in\mathbb{N}$ and $\Re(\beta)>0,\Re(\mu)>0$.\\

\textbf{Case XI.} Replacing $m$ by $2m$ and after that setting $\ell=k,a=1+k, b=1-2m-\gamma, \lambda=\mu$ and $s=\gamma$ in equation \eqref{eq(7.1)} and using Dixon's truncated summation theorem \eqref{eq(2.22)}, we obtain
\begin{eqnarray}\label{eq(7.10)}
&&\mathcal{M}\left\{e^{-\mu t}{_2F_2} \left[\begin{array}{r} -2m,1+k;\\ -2m-k, 1-2m-\gamma;\end{array}\mu t\right]_{2m};\gamma\right\}\nonumber\\
&&=\frac{\Gamma(\gamma)}{\mu^{\gamma}}\frac{(1+k)_{m}(\gamma)_{m}2^{2m}\left(\frac{1}{2}\right)_m(1+k+\gamma)_{2m}}{(1+k)_{2m}(\gamma)_{2m}(1+k+\gamma)_{m}},
\end{eqnarray}
where $1-2m-\gamma\in\mathbb{C}\setminus\mathbb{Z}_0^{-}; m,k\in\mathbb{N}$ and $\Re(\gamma)>0,\Re(\mu)>0$.\\

\textbf{Case XII.} Replacing $m$ by $2m$ and after that setting $\ell=k,a=1+k, b=-2m-k, \lambda=\mu$ and $s=1+k$ in equation \eqref{eq(7.1)} and using Dixon's truncated summation theorem \eqref{eq(2.23)}, we obtain
\begin{eqnarray}\label{eq(7.11)}
&&\mathcal{M}\left\{e^{-\mu t}{_2F_2} \left[\begin{array}{r} -2m,1+k;\\ -2m-k, -2m-k;\end{array}\mu t\right]_{2m};1+k\right\}\nonumber\\
&&=\frac{\Gamma(1+k)}{\mu^{1+k}}\frac{(1+k)_{m}(1+k)_{m}2^{2m}\left(\frac{1}{2}\right)_m(2+2k)_{2m}}{(1+k)_{2m}(1+k)_{2m}(2+2k)_{m}},
\end{eqnarray}
where $m,k\in\mathbb{N}$ and $\Re(\mu)>0$.\\

\textbf{Case XIII.} Replacing $m$ by $2m+1$ and after that setting $\ell=k,a=\gamma, b=-2m-\gamma, \lambda=\mu$ and $s=1+k$ in equation \eqref{eq(7.1)} and using Dixon's truncated summation theorem \eqref{eq(2.24)}, we obtain
\begin{eqnarray}\label{eq(7.16)}
\mathcal{M}\left\{e^{-\mu t}{_2F_2} \left[\begin{array}{r} -2m-1, \gamma;\\ -2m-1-k, -2m-\gamma;\end{array}\mu t\right]_{2m+1};1+k\right\}=0,
\end{eqnarray}
where $\gamma,-2m-\gamma\in\mathbb{C}\setminus\mathbb{Z}_0^{-}; m,k\in\mathbb{N}$ and $\Re(\mu)>0$.\\

\textbf{Case XIV.} Replacing $m$ by $2m+1$ and after that setting $\ell=k,a=1+k, b=-2m-k-1, \lambda=\mu$ and $s=1+k$ in equation \eqref{eq(7.1)} and using Dixon's truncated summation theorem \eqref{eq(2.25)}, we obtain
\begin{eqnarray}\label{eq(7.17)}
\mathcal{M}\left\{e^{-\mu t}{_2F_2} \left[\begin{array}{r} -2m-1,1+k;\\ -2m-1-k, -2m-1-k;\end{array}\mu t\right]_{2m+1};1+k\right\}=0,
\end{eqnarray}
where $m,k\in\mathbb{N}$ and $\Re(\mu)>0$.\\

\textbf{Case XV.} Replacing $m$ by $2m$ and after that setting $\ell=k,a=\alpha, b=1+\alpha+2m, \lambda=\mu$ and $s=1+\alpha+2m+k$ in equation \eqref{eq(7.1)} and using Dixon's truncated summation theorem \eqref{eq(2.20c)}, we obtain
\begin{eqnarray}\label{eq(7.12)}
&&\mathcal{M}\left\{e^{-\mu t}{_2F_2} \left[\begin{array}{r} -2m,\alpha;\\ -2m-k, 1+\alpha+2m;\end{array}\mu t\right]_{2m};1+\alpha+2m+k\right\}\nonumber\\
&&=\frac{\Gamma(1+\alpha+k+2m)}{\mu^{1+\alpha+k+2m}}\frac{(1+\alpha)_{2m}\left(1+\frac{\alpha}{2}+k\right)_{2m}}{\left(1+\frac{\alpha}{2}\right)_{2m}(1+k)_{2m}},
\end{eqnarray}
where $\alpha,1+\alpha+2m\in\mathbb{C}\setminus\mathbb{Z}_0^{-}; m,k\in\mathbb{N}$ and $\Re(1+\alpha+2m+k)>0,\Re(\mu)>0$.\\

\textbf{Case XVI.} Replacing $m$ by $2m+1$ and after that setting $\ell=k,a=\alpha, b=2+\alpha+2m, \lambda=\mu$ and $s=2+\alpha+2m+k$ in equation \eqref{eq(7.1)} and using Dixon's truncated summation theorem \eqref{eq(2.20d)}, we obtain
\begin{eqnarray}\label{eq(7.13)}
&&\mathcal{M}\left\{e^{-\mu t}{_2F_2} \left[\begin{array}{r} -2m-1,\alpha;\\ -2m-k-1, 2+\alpha+2m;\end{array}\mu t\right]_{2m+1};2+\alpha+2m+k\right\}\nonumber\\
&&=\frac{\Gamma(2+\alpha+2m+k)}{\mu^{2+\alpha+2m+k}}\frac{(1+\alpha)(2+2k+\alpha+4m)(2+\alpha)_{2m}\left(1+\frac{\alpha}{2}+k\right)_{2m}}{(2+\alpha)(1+2m+k)\left(2+\frac{\alpha}{2}\right)_{2m}(1+k)_{2m}},\nonumber\\
\end{eqnarray}
where $\alpha,2+\alpha+2m\in\mathbb{C}\setminus\mathbb{Z}_0^{-}; m,k\in\mathbb{N}$ and $\Re(2+\alpha+2m+k)>0,\Re(\mu)>0$.\\

\textbf{Remark.} In the next communication \cite{Qureshi2}, we shall obtain the Mellin transform of the product of exponential function and infinite Goursat series ${_2F_2} \left[\begin{array}{r} -m,\alpha;\\ -m-\ell, \beta;\end{array}\lambda t\right]$.

\section*{Concluding remarks}
In previous sections, we have derived some summation theorems for Clausen's terminating and truncated hypergeometric series ${_3F_2}$ when one numerator and one denominator parameters are negative integers. In the sequel of this paper, we have derived some summation formulae for Gauss' hypergeometric series ${_2F_1}$, Clausen hypergeometric series ${_3F_2}$ and have discussed their applications (see for example \cite{Qureshi2,Qureshi3}). It is expected that these summation formulae will be of wide interest and will help to advance research in the field of special functions.\\
We conclude our present investigation by observing that several hypergeometric summation theorems can be derived from a known summation theorem in an analogous manner.}

\end{document}